\documentclass[11pt,a4paper,reqno]{article}
\usepackage{amsmath,amsthm}

\def\br#1{\left(#1\right)}
\def\brs#1{\left\{#1\right\}}

\def\td{{\delta}}
\def\tdim#1#2{\Delta\left(#1,#2\right)}
\def\tdimu#1#2{\Delta^*\left(#1,#2\right)}

\newcommand{\supp}{\operatorname{supp}}

\newcommand{\loc}{\operatorname{loc}}

\usepackage{amsfonts}         
\usepackage{amsmath}
\usepackage{amssymb}

\newcommand{\R}{\mathbb{R}}
\newcommand{\C}{\mathbb{C}}
\newcommand{\D}{\mathbb{D}}

\newcommand{\onto}{\mathop{\longrightarrow}\limits}
\newcommand{\diam}{\operatorname{diam}}
\renewcommand{\Re}{\operatorname{Re}}

\swapnumbers
\theoremstyle{plain}
\newtheorem{thm}[equation]{Theorem}
\newtheorem{lem}[equation]{Lemma}
\newtheorem{prop}[equation]{Proposition}
\newtheorem{cor}[equation]{Corollary}

\theoremstyle{definition}

\newtheorem{conj}[equation]{Conjecture}

\theoremstyle{remark}
\newtheorem{rem}[equation]{Remark}

\numberwithin{equation}{section}

\title{ \textbf{\scshape Quasisymmetric distortion spectrum}}
\author{\textsc{Istv\'an Prause and Stanislav Smirnov}}
\date{}

\begin{document}

\maketitle

\footnotetext{\textit{Date:} October 25, 2009}
\footnotetext{
I.P. gratefully acknowledges funding from the Marie Curie Research training network
CODY of the European Commission.
This research was supported by projects 118634 and 1134757 of the Academy of Finland,
by the Swiss NSF and by the European Research Council AG CONFRA.}

\begin{abstract}
We give improved bounds for the distortion of the Hausdorff dimension under quasisymmetric maps 
in terms of the dilatation of their quasiconformal extension. 
The sharpness of the estimates remains an open question and 
is shown to be closely related to the fine structure of harmonic measure.
\end{abstract}

\pagestyle{plain}

\section{Introduction}

A homeomorphism $\phi \colon \Omega \to \Omega'$ between planar domains is called \emph{$k$-quasiconformal} if it 
belongs to the Sobolev class $W^{1,2}_{\loc}(\Omega)$ and satisfies the Beltrami equation 
\begin{equation}
\label{eq:beltrami}
\bar \partial \phi(z) = \mu(z) \partial \phi(z) \quad \mbox{a.e. } z \in \Omega,
\end{equation}
with a measurable coefficient $\mu$, $\| \mu \|_\infty \leq k <1$. An equivalent definition
says that infinitesimal circles are mapped to ellipses of eccentricity  bounded by 
the maximal dilatation $K$, with two constants related by $K=(1+k)/(1-k)$.

Quasisymmetric maps have been introduced by Beurling and Ahlfors \cite{beurling-ahlfors} as boundary correspondence under quasiconformal self-maps of the half-plane. An increasing function $g \colon \R \to \R$ is \emph{quasisymmetric} if
\[ \frac{1}{\rho} \leq \frac{g(x+y)-g(x)}{g(x)-g(x-y)} \leq \rho,
\]
for some constant $\rho \ge 1$ and for all $x,y \in \R$. We will call this the $\rho$\emph{-definition of quasisymmetry}.

Any quasisymmetric map admits quasiconformal extensions to the plane and thus we call a mapping 
of the real line \emph{$k$-quasisymmetric} if it can be extended to a $k$-quasiconformal map.
Without loss in dilatation constant $k$ we may 
require the extension to satisfy the reflection symmetry $\phi( z)=\overline{\phi(\bar z)}$. 
In what follows, we assume that every quasisymmetric map is endowed with some symmetric quasiconformal extension.

The two definitions are quantitatively equivalent:
a quasisymmetric map in the $\rho$-definition sense is 
$k(\rho)$-quasisymmetric with $k(\rho) \leq 1-1/\rho$, cf. \cite{lehtinen}.
Conversely, a $k$-quasisymmetric map satisfies the $\rho$-definition with  $\rho \leq 1/16 e^{\pi K}$.
In the latter direction even the best possible function is known,
which is given by a special function related to the hyperbolic metric of the three-punctured sphere \cite{lehto-virtanen}.
However, there is no exact correspondence between these two ways of quantifying quasisymmetry 
and therefore working with different definitions naturally leads to complementary results. 
In order to somewhat simplify the various expressions involved we make the convention to use quasiconformal dilatations $k\in[0,1[$ and $K \in[1,\infty[$ simultaneously.

Quasisymmetric maps need not be absolutely continuous with respect to the Lebesgue measure \cite{beurling-ahlfors}, in fact, a set of positive length may be compressed to a set of arbitrary low positive Hausdorff dimension \cite{tukia}. 
In the present note we estimate the structure of their singular sets in terms of multifractal spectra. 

Astala \cite{astala} gave optimal bounds for the dimensional distortion under general quasiconformal maps. 
We build upon his work, also employing holomorphic motions and thermodynamic formalism. 
In the case of quasisymmetric maps, the motion will be symmetric and this allows us to exploit extra information. Similarly to the general quasiconformal estimate, where extremal distortion is described by a conformal automorphism of the disk (cf. the function $\Phi$ in \eqref{eq:Phi}), we find bounds in terms of a degree two Blaschke product.
The paper is a direct follow-up to our earlier work \cite{prause,smirnov}.

Our main result is the estimate below for compression and expansion under quasisymmetric maps. 
In order to describe it, we introduce some notation. 
Hausdorff dimension will be referred to as ``$\dim$'' and dimension
of the sets will be usually denoted by the letter ``$\td$.'' 
Define the following function of $\td \in [0,1]$ and $k \in [-1,1]$:
\begin{equation} 
\label{eq:tofk}
 \tdim{\td}{k}:= \frac{\td\,(1-k^2)}{(1+k \sqrt{1-\td})^2} = 1-B_{-\sqrt{1-\td}}(k).
\end{equation}

The second expression is written in terms of the Blaschke product defined in \eqref{eq:blaschke} below.
The definition can also be nicely rewritten in terms of the minimal ``dilatations'' $\ell$ and $L$ of quasisymmetric maps such that our results possibly allow to take sets of dimension $1$ to dimension $\td$. Namely, for $\td=1-\ell^2$, one has
\begin{equation*} 
\label{eq:tofk-rewrite}
 \tdim{\td}{k}=\frac{(1-k^2)(1-\ell^2)}{(1+k\ell)^2}= 1-\br{\frac{k+\ell}{1+k\ell}}^2=
\frac{4 K L}{(KL+1)^2}= 1-\br{\frac{KL-1}{KL+1}}^2.
\end{equation*}
Denote the inverse function in $\td$ by $\tdimu{\td}{k}$, namely set
\begin{equation} 
\label{eq:tofku}
 \tdimu{\tdim{\td}{k}}{k}=\td~~\text{and}~~
 \tdim{\tdimu{\td}{k}}{k}=\min\brs{\td,1-k^2}.
\end{equation}
This can also be rewritten as
\begin{equation} 
\label{eq:tofkurewrite}
 \tdimu{\td}{k}=\tdim{\td}{- \min\{k,\sqrt{1-\td}\}}.
\end{equation}

\begin{thm}
\label{thm:qscompression}
Let $\phi \colon \R \to \R$ be a $k$-quasisymmetric map with some $k \in [0,1[$.
Then given a set $E \subset \R$ with Hausdorff dimension $\dim E=\td$ one has
\[ 
  \tdim{\td}{k} \leq \dim \phi (E) \leq \tdimu{\td}{k}.
\]
\end{thm}

In particular, setting $\td=1$ gives the statement: $k$-quasisymmetric image of the Lebesgue measure has dimension at least $1-k^2$, which is also discussed in \cite{prause}. 
The novelty of Theorem \ref{thm:qscompression} is the extension to the full spectrum of dimensions $\td \in [0,1]$. For related results in terms of the $\rho$-definition we refer to \cite{hayman-hinkkanen,heurteaux}.
Section 2 contains the proof of Theorem \ref{thm:qscompression} while in Section 3 we discuss various multifractal spectra and connections to harmonic measure.

The question of optimality of our estimates remains open and is directly related to the existence of
$k$-quasicircles with dimension $1+k^2$, see \cite{smirnov}. 
Moreover, there are intricate connections to the multifractal structure of harmonic measure, which will be the subject of our future work \cite{aps}.

\section{Quasisymmetric compression}

Our main technical result is the following

\begin{prop}
\label{prop:dd}
Let $\phi \colon \C \to \C$, $0\mapsto0$, $1\mapsto1$ be a $k$-quasiconformal map symmetric with respect to the real line: $\phi(z)=\overline{\phi(\bar z)}$.  Consider a collection of disjoint disks centered on the real line: 
$B_i=B(z_i,r_i) \subset \D$, $z_i \in \R$. Then for every $\rho \in (k,1)$, there exist a constant $a=a(\rho)>0$, such that following implication holds.
Assume that $\sum (ar_i)^\td \ge 1$ for some $\td \in [0,1]$, then
\[\sum (a \diam \phi B(z_i,r_i))^{\tdim{\td}{k/\rho}} \ge 1,
\]
where $\tdim{\td}{k}$ is defined by \eqref{eq:tofk}.
\end{prop}

\begin{proof}
The Beltrami coefficient $\mu$ of the map $\phi$ is symmetric with respect to the real axis: $\mu(z)=\overline{\mu(\bar z)}$. Embed $\phi$ into a holomorphic motion in a standard way: 
set the Beltrami coefficient $\mu_\lambda := \lambda \cdot \frac{\mu}{k}$ and denote the solution preserving $0,1$ and infinity by $\phi_\lambda(z)$. By the uniqueness of the solution we recover $\phi=\phi_k$. Solutions inherit the symmetry of $\mu$ in the form 
\begin{equation}
\label{eq:sym}
\phi_\lambda(z)=\overline{\phi_{\bar \lambda}(\bar z)}. 
\end{equation}
In particular, for real $\lambda$, $\phi_\lambda(\R)=\R$. 
Another crucial property for us is the holomorphic dependence of $\phi_\lambda$ on $\lambda$, see \cite{ahlfors-bers}.

We are interested in how the disks $\{ B(z_i,r_i) \}$ evolve in this motion. In order to have uniform estimates from now on we restrict the motion to the smaller disk $\rho \D$, with $k<\rho<1$. 
There is a constant $1\leq C=C(\rho)< \infty$ such that 
$$|x-z| \leq |y-z|~~\Rightarrow
~~|\phi_{\lambda}(x)- \phi_{\lambda}(z)| \leq C\,|\phi_{\lambda}(y)- \phi_{\lambda}(z)|
~\text{for~any}~|\lambda| \leq \rho.$$ 
This is the \emph{quasisymmetry property} of quasiconformal maps, see e.g. \cite{lehto-virtanen}.
In particular, under the restricted holomorphic motion the disks $B\br{\phi_\lambda(z_i),\frac1C |\phi_\lambda(z_i+ r_i)-\phi_\lambda(z_i)|}$ stay disjoint and included into $B(0,C)$. 
With one more rescaling we have a holomorphic family of disjoint disks 
\[ B\left(\frac1C \phi_\lambda(z_i),a |\phi_\lambda(z_i+ r_i)-\phi_\lambda(z_i)|\right)
\]
inside the unit disk.
In this step we choose the constant $a=a(\rho)=1/C^2$ in the statement of the proposition. 
We will work with the ``complex radius'' $r_i(\lambda)=a\br{\phi_\lambda(z_i+ r_i)-\phi_\lambda(z_i)}$ 
of these disks, in particular $r_i(0)=a r_i$. 
By a standard procedure this configuration of disks generates Cantor sets in $\D$,
which we denote by $C_\lambda$.

Recall the {\em variational principle} from \cite{astala} 
for the {\em pressure} $P_\lambda$ (in this elementary setting it is
a straight-forward application of the Jensen's inequality):
\begin{equation}
\label{eq:varprinciple}
P_{\lambda}(d):=\log \left( \sum |r_i(\lambda)|^{d} \right) =\, \sup_p (I_p - d \Re \Lambda_p(\lambda)),
\end{equation}
where for the probability distribution $\{ p_i \}$ we denote the {\em entropy} by
$$I:=- \sum p_i \log p_i,$$ 
and the {\em ``complex Lyapunov exponent''} by
$$\Lambda_p(\lambda)= - \sum p_i \log r_i(\lambda).$$ 
We fix the principal branch of the logarithm, which makes $\Lambda_p(\lambda)$ a holomorphic function in $\lambda$ for any fixed $p$. 
Recall that by the Bowen's formula the Hausdorff dimension $\dim C_\lambda$ is the unique root 
$\delta$ of the equation
$$P_{\lambda}(d)=0,$$
and so by the variational principle 
\begin{align}
\dim C_\lambda\le \delta~~&\Leftrightarrow~~P_\lambda(\delta)\ge0
~~\Leftrightarrow~~ \forall\,p,~I_p\ge \delta \Re\Lambda_p(\lambda)~,\label{eq:dimless}\\
\dim C_\lambda\ge \delta~~&\Leftrightarrow~~P_\lambda(\delta)\le0
~~\Leftrightarrow~~ \exists\,p,~I_p\le \delta \Re\Lambda_p(\lambda)~.\label{eq:dimmore}
\end{align}

Our task in terms of the pressure function is to show the following implication:
\[ P_0(\td) \ge 0 \quad \Longrightarrow \quad P_k(\tdim{\td}{k/\rho}) \ge 0.
\]
To this end, let us ``freeze'' $p$ at its value,
which maximizes  $P_0(\td)$ in the variational principle. 
Namely, by Jensen's inequality, set $p_i=r_i^\td/\sum r_i^\td$. 
Define the following holomorphic function
\begin{equation}
\label{eq:Phi} 
\Phi(\lambda)=1- \frac{I_p}{\Lambda_p(\lambda)}.
\end{equation}
In view of \eqref{eq:dimless} and the obvious $\dim C_\lambda \leq 2$, 
we have $I_p -2 \Re \Lambda_p(\lambda) \leq 0$, or equivalently $\Phi$ maps into the unit disk, 
$$\Phi \colon \rho \D \to \D.$$ 
Due to the symmetry \eqref{eq:sym} we have 
$$\Phi(\lambda)=\overline{\Phi(\bar \lambda)}.$$ 
Moreover, for real $\lambda$, all the disks are centered on the real line and hence $\dim C_\lambda \leq 1$, so by \eqref{eq:dimless} we have
$$\Phi(\lambda) \ge 0~\text{for}~\lambda \in \R.$$
Finally, by the choice of $p$ and our assumption $P_0(\td) \ge 0$ we have 
$$\Phi(0) \leq 1-\td.$$

In the next lemma we analyze the extremal problem described in the last paragraph and show that
\[
 \Phi(k) \leq B_{-\sqrt{1-\td}}(k/\rho) = 1- \tdim{\td}{k/\rho},
\]
therefore $I_p / \Lambda_p(k) \ge \tdim{\td}{k/\rho}$. Once more referring to \eqref{eq:varprinciple} we conclude the required estimate ${P_k(\tdim{\td}{k/\rho} \ge 0}$, thus proving the Proposition.
\end{proof}

\begin{lem}
Let $h \colon \D \to \D$, $h(z)=\overline{h(\bar z)}$ be a holomorphic map, sending the interval $(-1,1)$ into $[0,1)$.
Suppose that $h(0) \leq l^2$ for some $l \ge 0$. Then for any $k \in [0,1)$ 
\[ h(k) \leq \br{\frac{k+l}{1+k l}}^2.
\]
The extremal map is given by the degree two Blaschke product $B_{-l}$ with an order two zero at $-l$:
\begin{equation}
\label{eq:blaschke}
B_{-l}(z)=\left( \frac{z+l}{1+l z} \right)^2.
\end{equation}
\end{lem}

\begin{proof}
If $l=0$ then our assumptions force $h(0)=0$ and $h'(0)=0$. Therefore Schwarz lemma applied to $h(z)/z$ implies $h(k) \leq k^2$ and the lemma. Alternatively, we may reduce the $l=0$ case to $l>0$ by considering the limit $l \to 0$. 

From now on we assume $l >0$. We will use a three-point version of Schwarz-Pick lemma by Beardon and Minda \cite{beardon-minda} in this case. Let us first briefly recall their argument. For $z,w \in \D$, set
\begin{equation}\label{eq:hdef}
 [z,w]=\frac{z-w}{1-\bar w z}\in\D,\quad \quad h^*(z,w)=\frac{[hz,hw]}{[z,w]}\in \overline{\D}.
\end{equation}
The first quantity relates to the hyperbolic distance in the following way:
\begin{equation}\label{eq:hyp} 
\big|\,[z,w]\,\big| = \tanh\br{\frac{d(z,w)}2}.
\end{equation}
For a fixed $w$, the function $h^*$ is holomorphic in $z$. 
Since $h$ is not a conformal automorphism of $\D$, by the standard Schwarz-Pick lemma, 
$h^*$ maps holomorphically \emph{into} the unit disk.
Yet another application of the Schwarz-Pick lemma gives (see \cite[Theorem 3.1]{beardon-minda})
the following inequality for the hyperbolic distance:
\begin{equation}
\label{eq:3point} 
d\br{h^*(z,v),h^*(w,v)} \leq d\br{z,w}. 
\end{equation}
We will use this three-point Schwarz lemma for $z=k$, $v=0$, and $w=-l$.

If $h^*(k,0) \le 0$, then $h(k)-h(0)\le0$ by \eqref{eq:hdef}, and the Lemma follows from
$$h(k)\le h(0)\le l^2\le\br{\frac{k+l}{1+k l}}^2.$$
So from now on we assume that $h^*(k,0) \ge 0$. 

Now if $h^*(-l,0)\le0$, we can write
$$d(0,h^*(k,0)) \leq d(h^*(-l,0),h^*(k,0)) \leq d(-l,k) = d(0,k)+d(0,l),
$$
implying \eqref{eq:hstar} below.
If, on the contrary, $h^*(-l,0)\ge0$, then by \eqref{eq:hdef} we have $h(-l)\le h(0)$.
Combining this with our assumptions we arrive at
$$0\le h(-l)\le h(0)\le l^2,$$
and therefore
$[h(-l),h(0)]\le l^2$.
We conclude that
\begin{equation}
\label{eq:l}
 h^*(-l,0)=\frac{[h(-l),h(0)]}{[-l,0]}\leq \frac{l^2}{l}=l.
\end{equation}
Combining \eqref{eq:3point} and \eqref{eq:l} we may write
\begin{align*}
d(0,h^*(k,0)) &\leq d(0,h^*(-l,0)) + d(h^*(-l,0),h^*(k,0))\\ 
&\leq d(0,l)+d(-l,k) = d(0,k)+2d(0,l),
\end{align*} 
concluding that
\begin{equation}\label{eq:hstar} 
d(0,h^*(k,0))\leq d(0,k)+2d(0,l)=d\br{0,\frac{k+2l+kl^2}{1+2kl+l^2}}. 
\end{equation}
Therefore 
\begin{equation}\label{eq:hk} 
\frac{[h(k),h(0)]}{[k,0]}=h^*(k,0)\leq \frac{k+2l+kl^2}{1+2kl+l^2}, 
\end{equation}
which together with \eqref{eq:hyp} gives an upper estimate for $d(h(0),h(k))$ 
in terms of only $k$ and $l$.
A direct calculation shows that for $h:=B_{-l}$
we have an equality in \eqref{eq:hk},
and thus
$$d(h(0),h(k)) \leq d(B_{-l}(0),B_{-l}(k)).$$
Then using $0\le B_{-l}(0)\le B_{-l}(k)$ and
$h(0) \leq l^2 = B_{-l}(0)$ we deduce
\begin{align*}
 d(0,h(k))&= d(0,h(0)) + d(h(0),h(k)) \\
&\leq d(0,B_{-l}(0)) + d(B_{-l}(0),B_{-l}(k)) = d(0,B_{-l}(k)).
\end{align*}
We have shown that $h(k) \leq B_{-l}(k)$, as required.
\end{proof}

\noindent
Theorem \ref{thm:qscompression} easily follows from Proposition \ref{prop:dd}:

\begin{proof}[Proof of Theorem \ref{thm:qscompression}]
First of all, it is sufficient to consider the lower estimate as the upper bound follows from considering the inverse map and using (\ref{eq:tofku}). 
Proposition \ref{prop:dd} establishes the required compression relation for disjoint packings. 
The normalization assumption in Proposition \ref{prop:dd} does not influence dimension estimates 
and we may also assume $E \subset [-1,1]$. It is a routine application of the covering theorems 
(e.g. the $5r$-covering lemma) 
to pass from coverings to packings and hence we conclude that if $E$ is of infinite $\td$-dimensional Hausdorff measure then $\dim \phi(E) \ge \tdim{\td}{k/\rho}$. Sending $\rho \to 1$ we find that $\dim \phi(E) \ge \tdim{\td}{k}$. Finally, a limiting argument in $\td$ shows that the
conclusion holds also if we only assume $\dim E=\td$. For more details on the covering argument, we refer the reader to \cite{prause}.
\end{proof}

\section{Multifractal spectra}

Multifractal analysis of {harmonic measure} provides a suitable
framework for discussing compression and expansion phenomena of conformal maps. 
In the rest of the paper, we discuss analogous multifractal spectra for {quasisymmetric} maps and 
show how multifractality of harmonic measure is reflected in the singularity of the welding.

The procedure of \emph{conformal welding} gives rise to a correspondence between Jordan curves and homeomorphisms of the unit circle in the following manner. 
Given a closed Jordan curve $\Gamma$, let $g_+ \colon \D \to \Omega$ and $g_- \colon  \hat \C \setminus \overline{\D} \to \Omega^*$ be conformal maps onto the bounded and unbounded complementary components of $\Gamma$ respectively. 
Then the boundary correspondences induce a homeomorphism $\phi=g_-^{-1} \circ g_+ \colon \partial \D \to \partial \D$, 
and a homeomorphism arising this way is called a conformal welding. 
Given a homeomorphism $\phi$ of the unit circle, major open problems are
to understand whether it is a conformal welding, to find the corresponding
curve $\Gamma$, and to determine whether it is unique.

A powerful tool for solving the welding problem is the Beltrami equation \eqref{eq:beltrami}.
The situation is well-understood in the uniformly elliptic setting: quasisymmetric maps are conformal weldings and in this case $\Gamma$ is a quasicircle, see e.g.~\cite{aim,lehto-virtanen}.
More quantitatively, the welding is $K^2$-quasisymmetric if and only if the conformal map $g_+$ admits a $K^2$-quasiconformal extension. Note that by \cite{smirnov} the latter is equivalent to $\Gamma$ being a $K$-quasicircle.
To summarize, we have the following exact correspondence:
\[\text{welding $\phi$ is $K^2$-quasisymmetric} \iff \ \text{$\Gamma$ is a $K$-quasicircle}.
\]

\subsection*{Quasisymmetric spectra}

In this section we rephrase our results using the language of \emph{multifractal analysis}. 
Let $\phi \colon \D \to \D$ be a $K$-quasiconformal self-map of the unit disk. 
We are going to analyze the multifractal structure of the push-forward $\mu=\phi_*(m)$ of the 
normalized Lebesgue measure $m$ on $\partial \D$.

\paragraph{Box dimension spectrum.}
This spectrum describes the size of the set where the measure scales with a fixed exponent $\alpha>0$. More precisely, we define
\[ f_\mu(\alpha)= \lim_{\epsilon \to 0} \limsup_{{r} \to 0} \frac{\log N({r},\alpha,\epsilon)}{|\log {r}|},
\]
where $N({r},\alpha,\epsilon)$ is a maximum number of disjoint disks $B_n=B(z_n,{r})$
with centers $z_n \in \supp \mu$ and ${r}^{\alpha+\epsilon} \leq \mu (B_n) \leq {r}^{\alpha-\epsilon}$.
We also define the spectrum for the class of $K$-quasiconformal self-maps of the disk as
\[ F_{K\mbox{-qs}}(\alpha)= \sup \{ f_{\phi_*(m)}(\alpha) \ | \ \phi \colon \D 
\onto^{\mbox{\tiny onto}} \D \ K\mbox{-quasiconformal} \}.
\]
Observe that by H\"older continuity this spectrum is equal to $-\infty$ outside the interval $[1/K,K]$,
if we use the usual convention $\dim \emptyset= -\infty$.

\paragraph{Integral means spectrum.}
The \emph{integral means spectrum} of $\phi$ and its universal counterpart are defined respectively as 
\begin{equation*}
\beta_\phi(t)=\inf \left\{ \beta \colon \int_0^{2\pi} \left(\frac{1-|\phi(re^{i\theta}|}{1-r} \right)^t d \theta =O( (1-r)^{-\beta}) \right\}, \quad t \in \R,
\end{equation*}
\[ B_{K\mbox{-qs}}(t)=\sup \{ \beta_\phi(t) \ | \ \phi \colon \D \onto^{\mbox{\tiny onto}} \D \ K\mbox{-quasiconformal} \}.
\]
The definitions above are  motivated by the corresponding notions for harmonic measure and  univalent maps, see Makarov's \cite{makarov}. 
The integral means spectrum $\beta_g(t)$ for a conformal map $g \colon \D \to \Omega$, for instance, 
is defined in the same manner except the difference quotient is replaced by $|g'(re^{i \theta})|$. 

\begin{thm}
\label{thm:qsspectra}
The following upper bounds hold true:
\begin{align*}
 F_{K\mbox{-qs}}(\alpha) &\leq - \frac{4K}{(K-1)^2} (\sqrt{\alpha}- \sqrt{K})(\sqrt{\alpha}- 1/\sqrt{K})\\
\text{for}&~\frac1K \leq \alpha \leq 1-k^2~~\text{and}~~\frac{1}{1-k^2} \leq \alpha \leq K.
\end{align*}
In the range $1-k^2 \leq \alpha \leq 1/(1-k^2)$ we have the trivial bound
\[
 F_{K\mbox{-qs}}(\alpha) \leq \min \{ \alpha,1 \}.
 \]
The $\beta$-spectrum satisfies the following estimate
\[ B_{K\mbox{-qs}}(t) \leq \max \left\{0, \frac{t(t-1)}{t+\frac{4K}{(K-1)^2}} \right\} \quad \mbox{for } 
-\frac{2}{K-1} \leq t \leq  \frac{2K}{K-1}.
\]
At the endpoints a phase transition occurs and the spectrum becomes linear,
\[ B_{K\mbox{-qs}}(t) = \left\{
 \begin{array}{ll}
 -(K-1)t-1 & \mbox{for $t \leq -2/(K-1)$}, \\
(1-1/K)t-1 &  \mbox{for $t \ge 2K/(K-1)$}.
 \end{array} \right.
\]
\end{thm}

\begin{rem}
The linear part was already established by Bishop \cite{bishop}. He also discussed  
possible general values of $B_{K\mbox{-qs}}(t)$, see his Questions 5.6 and 5.7. 
The proposed function has vanishing left-derivative at $0$ which is not in agreement with 
the fact that sets of full dimension may be compressed under quasisymmetric map. 
Instead, we conjecture that the bounds in Theorem \ref{thm:qsspectra} are
in fact optimal.
\end{rem}

\begin{proof}
The estimates are rather direct consequences of Proposition \ref{prop:dd}. 
We only sketch the calculations for the $f$-spectrum, the $\beta$-spectrum is obtained via a Legendre transform. 
First, let us record the symmetry due to the invariance under inverse mappings:
\[  F_{K\mbox{-qs}}(\alpha) = \alpha  F_{K\mbox{-qs}}\left(\frac{1}{\alpha} \right),
\]
or in terms of the $\beta$-spectrum:
\[  B_{K\mbox{-qs}}(  B_{K\mbox{-qs}}(t)-t+1)=  B_{K\mbox{-qs}}(t).
 \]
Thus, it suffices to consider the compression case $\alpha \leq 1$.
Formally, $ F_{K\mbox{-qs}}(\alpha) \leq \tdim{\td}{k}$ when $\td$
is chosen so that $\tdim{\td}{k}=\alpha\td$.
Indeed, the compression bound of Proposition \ref{prop:dd}, 
applied to the collection of disks in the definition of
the box dimension spectrum and to the inverse map $\phi^{-1}$
provides for all $\td$ the implication
\[  
F_{K\mbox{-qs}}(\alpha) \geq \alpha\td ~~\Rightarrow~~
F_{K\mbox{-qs}}(\alpha) \geq \tdim{\td}{k},
\]
so if the first inequality holds for some $\td$ with $\alpha\td<\tdim{\td}{k}$,
we recursively deduce that $F_{K\mbox{-qs}}(\alpha)$ is infinite,
arriving to a contradiction. 
Therefore
\[  F_{K\mbox{-qs}}(\alpha) \leq \sup \{ \tdim{\td}{k} : \alpha \td \ge \tdim{\td}{k}\}.
\]
Solving $\tdim{\td}{k}= \alpha \td$ leads to the expression in the statement of the theorem.
\end{proof}

\subsection*{Spectrum of quasidisks}

In this section we relate the quasisymmetric spectrum to the spectra of conformal maps with quasiconformal extensions. 
Our main tool is the decomposition of a quasiconformal map into 
\emph{symmetric} and \emph{antisymmetric} parts from \cite{smirnov}, which we recall below.

It will be more convenient to consider symmetry with respect to the real line, so let us assume
that $g \colon \C \to \C$ is a $K^2$-quasiconformal map which is conformal in the upper half-plane 
$\C_+$.  Then $g$ can be written \cite[Theorem 2]{smirnov} as a superposition $g=\psi \circ \phi$, 
where $\phi$ and $\psi$ are global $K$-quasiconformal mappings and $\phi$ is symmetric while $\psi$ is antisymmetric with respect to the real line. 
The antisymmetry means that the Beltrami coefficient satisfies
\[ \mu_\psi(z)= - \overline{\mu_\psi(\bar z)}.
\]

\paragraph{Universal spectrum.}
Let us recall some definitions from \cite{makarov} for the relevant classes.
The \emph{universal integral means spectrum}  (for bounded univalent functions) is defined by 
\[ B(t)=\sup \{ \beta_g(t) \ | \ h \colon \D \to \Omega\subset\C \mbox{ is a bounded univalent map}  \}.
\]
Consider now the class of \emph{univalent maps with $K$-quasiconformal extension} and denote
its integral means spectrum by $B_K(t)$: 
\[ B_{K}(t)=\sup \{ \beta_g(t) \ | \ h \colon \D \to \Omega \subset\C
\mbox{ is univalent with a $K$-quasiconformal extension to $\C$} \}.
\]
As we already pointed out, a univalent map admits a $K^2$-quasiconformal extension if and only if the image domain is a ${K}$-quasidisk.

The \emph{universal spectrum conjecture} states \cite{kraetzer} that 
\begin{equation*}
 B(t)=\frac{t^2}{4} \quad \mbox{for } |t| \leq 2.
\end{equation*}
A somewhat stronger variant says that
\begin{equation}
\label{eq:qct2over4}
 B_K(t)=\frac{k^2 t^2}{4} \quad \mbox{for } |t| \leq \frac{2}{k}.
\end{equation}
We refer to \cite{beliaev-smirnov,jones} for further discussion of these conjectures.

\begin{thm}
\label{thm:quasidiskspectrum}
We have the upper bound
\begin{align}
\label{eq:qdiskspectrum1}
B_{K^2}(t) &\leq \frac{k^2 t^2}{(1+k^2)^2}=\frac14 \left( \frac{K^2-1}{K^2+1} \right)^2 t^2 \quad \mbox{for } 1+k^2 \leq t \leq \frac{1+k^2}{k},\\
\label{eq:qdiskspectrum2}
B_{K^2}(t) &= \frac{K^2-1}{K^2+1} t-1 \quad \mbox{for } t \ge \frac{1+k^2}{k}=2 \frac{K^2+1}{K^2-1}.
\end{align}
In other words, the conjectural upper bound in \eqref{eq:qct2over4} holds from the point $1+k^2$
onwards for the spectrum of $K$-quasidisks.
\end{thm}

\begin{proof}
Let $g \colon \D \to \Omega$ be a conformal map with $K^2$-quasiconformal extension.
The theorem follows from the following dimension distortion bounds: if 
$E \subset \partial \D$, $\dim E=\delta$ then
\begin{equation}
\label{eq:expansion}
\dim g E \leq \frac{(1+k^2)\delta}{1+k^2-2k \sqrt{1-\delta}},
\end{equation}
provided $\delta \leq 1-k^2$. Otherwise we have $\dim g E \leq \dim \partial \Omega \leq 1+k^2$ from \cite{smirnov}.
In fact, the integral means bound \eqref{eq:qdiskspectrum1} corresponds to the dimension expansion bound \eqref{eq:expansion} for Gibbs measures via Legendre transform \cite{makarov}. In the range 
$\delta \leq 1-k^2$ we prove the stronger fact that the expansion bound holds for arbitrary sets.

In order to prove \eqref{eq:expansion} we transfer the setting from the unit disk to the half-plane.
Let us adjust our notation: $K^2$-quasiconformalmap $g \colon \C \to \C$ is  conformal in $\C_+$ 
and set $E \subset \R$ has dimension $\dim E=\delta$.
Apply now the decomposition into the symmetric and antisymmetric parts, $g=\psi \circ \phi$, as above.
We will use the expansion estimates separately for both $\phi$ and $\psi$. 
Expansion by the map $\phi$ is estimated by Theorem \ref{thm:qscompression}:
\begin{equation*}
\dim \phi(E) \leq \tdimu{\delta}{k}=:\Delta^*.
\end{equation*}
For the other map $\psi$ we make use of the improvement distortion estimates from \cite{smirnov}. 
The precise bound is given by (see also \cite[Theorem 13.3.6]{aim})
\begin{equation}
\label{eq:expasym}
\dim \psi (\phi (E)) \leq \frac{(1+k^2)\Delta^*}{1-k^2+k^2 \Delta^*}.
\end{equation}
Finally, substituting \eqref{eq:tofkurewrite} for $\Delta^*$ yields \eqref{eq:expansion}.
The equality in \eqref{eq:qdiskspectrum2} follows from considering a domain whose boundary has an angle-type singularity, such as in \cite{becker-pommerenke}.
\end{proof}

Let us point out the end-point version of the previous theorem in an integrability form.
Previously, it was known that in the class of univalent functions which admit a
$K$-quasiconformal extension the H\"older exponent improves from $1/K$ to $1-k=2/(K+1)$, 
see \cite{becker-pommerenke,pommerenke}.
We show that this improvement holds true on the level of integrability of the derivative.
The refinement is to be compared to the exponent $2K/(K-1)$ for general quasiconformal mappings \cite{astala}.

\begin{cor}
If $\phi \colon \D \to \C$ is a conformal map with $K$-quasiconformal extension then
\[ \phi' \in L^p(\D) \quad \mbox{for all } \quad 2 \leq p < \frac{2(K+1)}{K-1}.
\]
The upper bound for the exponent is the best possible.
\end{cor}

\begin{proof}
This is equivalent to the statement $B_K(2/k)=1$ of Theorem \ref{thm:quasidiskspectrum}.
\end{proof}

\paragraph{Lower bounds.}
Astala asked in \cite{astala} whether $1+k^2$ is the correct bound on the dimension of quasicircles.
The work \cite{smirnov} confirmed the upper estimate, leaving only the question of sharpness open.
We formulate this as Astala's conjecture, see also \cite[Conjecture 13.3.2]{aim} for a discussion.
\begin{conj}[Astala's conjecture]
\label{conj:astala}
For every $0<k<1$, there exists a $k$-quasicircle $\Gamma$ with Hausdorff dimension
\[ \dim \Gamma =1+k^2.
\]
\end{conj}

\begin{rem}
Our proofs, as well as \cite{smirnov}, rely crucially on the Schwarz lemma. The rigidity property of the Schwarz lemma has the following consequence. If there exists a $k$-quasicircle with dimension 
$1+k^2$ for \emph{some} $0<k<1$ then in fact there exist quasicircles with the same property for 
\emph{all} $k$. Moreover, in the dual direction this implies that Theorem \ref{thm:qscompression} is sharp for all values of $t$ and $k$.
\end{rem}

A further connection between the quasisymmetric and conformal spectra is given by the following conditional theorem:

\begin{thm}
Astala's conjecture on quasicircles implies the conjectured lower bound in \eqref{eq:qct2over4} for negative $t$, that is
\begin{equation}
\label{eq:lowerbound}
B_K(t) \ge \frac{k^2 t^2}{4} \quad \mbox{for } -\frac{2}{k} \leq t \leq 0,
\end{equation}
for \emph{all} $K \ge 1$. In particular, it implies $B(t) \ge t^2/4$ for $t \in [-2,0]$.
\end{thm}

\begin{proof} 
As we remarked, Conjecture \ref{conj:astala} implies for any $\delta \in [0,1]$ and $K \ge 1$ the existence of a symmetric $K$-quasiconformal map $\phi$ which sends a set $E \subset \R$ of dimension $\delta \in [0,1]$ to a set of dimension $\Delta=\tdim{\delta}{k}$. 
Based on this quasisymmetric map we are going to produce a conformal map with strong contraction properties. One could produce such a map via the welding construction, here we use a related procedure.
Consider the inverse $\phi^{-1}$ and its Beltrami coefficient $\mu(z)$ in $\C_+$. 
Set the same coefficient $\mu(z)$ in $\C_+$, and extend it to $\C_-$ in an antisymmetric fashion, 
by the formula $-\overline{ \mu (\bar z)}$. 
The solution to this Beltrami equation is a $K$-quasiconformal antisymmetric map $\psi$, moreover the composition $g=\psi \circ \phi$ is a $K^2$-quasiconformal map
which is conformal in $\C_+$ by construction.

We can again apply the expansion bound $\eqref{eq:expasym}$ to the antisymmetric map $\psi$
and find that
\begin{equation}
\dim \psi (\phi (E)) \leq \frac{(1+k^2)\Delta}{1-k^2+k^2 \Delta}.
\end{equation}
Substituting \eqref{eq:tofk} in place of $\Delta$, we obtain the required contraction property
\[ \dim gE \leq \frac{(1+k^2) \delta}{1+k^2+2k \sqrt{1-\delta}} \quad \mbox{with } \quad \dim E =\delta.
\]
By conjugating with a M\"obius transformation we may transfer to the unit disk,  and find a map $g_* \colon \D \to \C$ satisfying the same bounds.
This exactly corresponds to \eqref{eq:lowerbound} via Makarov's formula \cite{makarov87}
\[ \dim g_* E \ge \frac{-t \dim E}{\beta_{g_*}(t)-t+1-\dim E}, \quad \mbox{for all } \quad t<0.
\]
\end{proof}

\begin{rem}
The restriction for negative values of $t$ in the previous theorem is not essential.
If we appropriately interpret $B(t)$ with complex $t$ (cf. \cite{binder}),
the conditional conclusion $B(t) \ge |t|^2/4$ holds for
complex values $|t| \leq 2$ as well. 
This requires a more precise understanding of the connection between 
harmonic measure and quasiconformal mappings, and is one of the subjects to be discussed in \cite{aps}.
\end{rem}

\noindent 
{\scshape Department of Mathematics and Statistics, University of Helsinki, Finland\\
P.O.~Box 68 FIN-00014 \\
E-mail: {\sf Istvan.Prause@helsinki.fi}}

\noindent 
{\scshape Section de Math\'ematiques, Universit\'e de Gen\`eve.
2-4 rue du Li\`evre, Case postale 64, CH-1211 Gen\`eve 4, Suisse\\
E-mail: {\sf Stanislav.Smirnov@unige.ch}}

\end{document}